\documentclass[12pt,leqno]{article}
\usepackage{amsmath, amssymb, amsfonts, amsthm,amscd}

\swapnumbers
\theoremstyle{plain}
\newtheorem{theorem}{Theorem}

\newtheorem{proposition}{Proposition}

\theoremstyle{definition}
\newtheorem{remark}[theorem]{Remark}

\numberwithin{equation}{theorem}
  \newcommand{\om}{\omega}
  \newcommand{\sig}{\sigma}

  \newcommand{\Hom}{\operatorname{Hom}}
 \newcommand{\reg}{\operatorname{reg}}
  \newcommand{\sing}{\operatorname{sing}}

\newcommand{\bfi}{\bf\itshape}
\newcommand{\subneq}{\underset{\neq}{\subset}}
\renewcommand{\ni}{\noindent}

 \newcommand{\ical}{\mathcal{I}}

 \newcommand{\ocal}{\mathcal{O}}

\begin{document}

  \title{Finiteness of the number of compatibly-split subvarieties}
  \author{Shrawan Kumar and Vikram B. Mehta\thanks{{\bf Acknowledgements.} We thank A. 
Knutson for the correspondences which led to this work.  
The first author 
was partially supported by a FRG grant from NSF.}}
  \maketitle

\section{Introduction}

Let $k$ be an algebraically closed field of characteristic $p>0$ and let 
$X$ be a scheme over $k$ (always assumed to be separated of finite type 
over $k$).  The following is the main theorem of this note and we give 
here its complete and self-contained proof. 

  \begin{theorem}
  Assume that $X$ is Frobenius split by a splitting $\sig\in$ \newline 
$\Hom_{\ocal_X} (F_*\ocal_X, \ocal_X)$, where $F$ is the absolute 
Frobenius morphism (cf. [BK, Section 1.1]).  Then, there are only finitely 
many closed subschemes of $X$ which are compatibly split (under $\sig$).
    \end{theorem}

  \section{Proof of Theorem 1.1}

We first prove the following proposition which is of independent interest.  
By a variety we mean a reduced 
but not necessarily irreducible scheme over $k$.

\begin{proposition} Let $X$ be a nonsingular irreducible variety which is Frobenius split by 
$\sig\in\Hom_{\ocal_X}(F_*\ocal_X, \ocal_X) \simeq H^0(X, F_*(\om_X^{1-p}))$, where $\om_X$ 
is the dualizing sheaf of $X$ (cf. [BK, Proposition 1.3.7]), and let $Y\subneq X$ be a 
compatibly-split closed subscheme of $X$.  Then,
  \[ Y \subset Z(\bar{\sig}),
  \] where $Z(\bar{\sig})$ denotes the set of zeroes of $\bar{\sig}$ and $\bar{\sig}$ is the 
section of $F_*(\om_X^{1-p})$ obtained from $\sig$ via the above identification.
  \end{proposition}

  \begin{proof} Since any irreducible component of a compatibly-split closed subscheme is 
compatibly split (cf. [BK, Proposition 1.2.1]), we can assume that $Y$ is irreducible.  
Assume, if possible, that
   $Y\cap (X\backslash Z(\bar{\sig})) \neq \emptyset .$ Then, $Y^{\reg} \cap 
(X\backslash Z(\bar{\sig}))\neq\emptyset ,$
 where $Y^{\reg}$ is the nonsingular locus of $Y$.

Take $y\in Y^{\reg}\cap (X\backslash Z(\bar{\sig}))$.  Choose a system of local parameters 
$\{ 
t_1,\dots ,$\\ $t_m, t_{m+1},\dots ,t_n\}$ at $y\in X$ such that $\{ t_1,\dots ,t_m\}$ is 
a system of local parameters at $y\in Y$ and 
$\langle t_{m+1},\dots ,t_n\rangle$ is the completion of the ideal of $Y$ in $X$ at $y$.   
  (This is possible since both $X$ and $Y$ are 
nonsingular at $y$.)  By assumption, $\bar{\sig}$ is a unit in the local ring $\ocal_{X,y}$.  
Moreover, $\sig$ induces a splitting $\hat{\sig}$ of the power series ring $k[[t_1,\dots 
,t_n]]$ compatibly splitting the ideal $\langle t_{m+1}, \dots , t_n\rangle$.  Now, since 
$\bar{\sig}$ does not vanish at $y$, $\hat{\sig}\bigl( (t_{1}\cdots t_n)^{p-1}\bigr)$ is a 
unit in the ring $k[[t_1,\dots ,t_n]]$.  In particular, $\hat{\sig}$ does not keep the ideal 
$\langle t_{m+1},\dots ,t_n\rangle$ stable.  This is a contradiction to the assumption.  
Hence, $Y\subset Z(\bar{\sig})$, proving the proposition.
  \end{proof}

\ni ({\bf 2.2}) {\bfi Proof of Theorem 1.1}.  By [BK, Proposition 1.2.1], we can assume 
without loss of generality that $X$ is irreducible.   
We prove Theorem 1.1 by induction on the dimension of $X$ . If dim $X = 0$, then the theorem 
is 
clear. So assume that dim $X = n$  and  the theorem is true for varieties of 
dimension $< n$.  Let $Y\subneq X$ be a compatibly-split irreducible closed subscheme.  
Then, 
either $Y\subset X^{\sing}$ (where $X^{\sing}$ is the singular locus of $X$) or $Y\cap 
X^{\reg}\neq\emptyset$.  In the latter case, by Proposition 2.1,
  \[ Y\cap X^{\reg} \subset Z(\bar{\sig^o}),
  \] where $Z(\bar{\sig^o})$ denotes the set of zeroes of the splitting 
$\bar{\sig^o}$ of the open subset $X^{\reg}$ of $X$ viewed as a section of 
$F_*\Bigl(\om^{1-p}_{X^{\reg}}\Bigr)$.  Thus, in this case,
  \[ Y\subset \overline{Z(\bar{\sig^o})} ,
  \] $\overline{Z(\bar{\sig^o})}$ being the closure of $Z(\bar{\sig^o})$ in $X$.  
Hence, in either case,
  \begin{equation} Y\subset \overline{Z(\bar{\sig^o})} \cup X^{\sing} .
  \end{equation}

  Considering the irreducible components, the same inclusion (1) holds for any 
compatibly-split closed subscheme $Y\subset X$ such that $Y\neq X$.

  Let $\{ Y_i\}_{i\in I}$ be the collection of all the distinct compatibly-split closed 
subschemes $Y_i\subneq X$ and let $Y := \overline{\bigcup_{i\in I} Y_i}$.  Since the ideal 
sheaf $\ical_Y = \bigcap_{i\in I} \ical_{Y_i}$ and each $\ical_{Y_i}$ is stable under the 
splitting $\sig$ of $X$, the closed subscheme $Y$ is compatibly split.  In particular, by 
(1), for each $i\in I$,
    \[ Y_i\subset Y\subset \overline{Z(\bar{\sig^o})} \cup X^{\sing}.
  \] Since $\dim \bigl( \overline{Z(\bar{\sig^o})}\cup X^{\sing}\bigr) < \dim X$; in 
particular, one has $\dim Y < \dim X$.  Thus, by the induction hypothesis (applying the 
theorem with
$X$ replaced by $Y$), $I$ is a finite set.  This completes the proof of the theorem.\\ 
\hspace*{5in} $\square$

\begin{remark}  
Karl Schwede has also obtained the above theorem via `$F$-purity' in a 
recent preprint [S].  
As pointed out by Schwede, when X is projective, the theorem 
also follows from [EH, Corollary 3.2] again via `$F$-purity'.
(See also, [Sh] for another proof via `tight closure'.)
\end{remark}

\medskip

\ni {\bf Addresses:}
\smallskip

\ni S.K.:  Department of Mathematics, University of North Carolina, Chapel Hill, NC 
27599-3250, USA\\ \quad (email:  shrawan$@$email.unc.edu)

\ni V.M.: School of Mathematics, Tata Institute of Fundamental Research,\\ Calaba, Mumbai 
400005, India \\ \quad (email: vikram$@$math.tifr.res.in)

  \end{document}